\documentstyle[amssymb,12pt]{article}
\title{Coding into $K$ by reasonable forcing}
\author{Ralf-Dieter Schindler\footnote{The author would like to thank
Itay Neeman,  
Philip Welch, and in particular Sy Friedman
for their interest and for their many
hints and comments.
John Steel even provided a crucial subclaim, and again
I do say thanks
for his intellectual support 
during my stay in Berkeley. I gratefully acknowledge
financial support from the Deutsche Forschungsgemeinschaft (DFG).}}
\newtheorem{defn}{Definition}[section]
\newtheorem{thm}[defn]{Theorem}

\newcommand{\forces}{\ \ |\!\!\!|\!- \ }

\begin{document}
\maketitle

\section{Introduction.}

The present paper was inspired by a talk Itay Neeman gave on his joint
work \cite{NZ1} and \cite{NZ2} with J. Zapletal. Assuming that,
vaguely, $AD^{L({\Bbb R})}$ holds (i.e., that $L({\Bbb R})$ is a model
of the Axiom of Determinacy) 
they can show that no
set of ordinals not already in $L({\Bbb R})$ can be 
coded into $L({\Bbb
R})$ by a set-sized reasonable forcing and moreover that
the theory of $L({\Bbb R})$ with
parameters for ordinals and reals from $V$ is
frozen with respect to
all generic extensions by set-sized reasonable forcings. 
"Reasonability" was
introduced by Foreman and Magidor in \cite{FM}. A notion of
forcing $P$ is called reasonable if for any infinite
ordinal $\alpha$ it is true that $[\alpha]^\omega \cap
V$ is stationary in $[\alpha]^\omega \cap V[G]$, for every $G$
being $P$-generic over $V$
(cf. \cite{FM} Definition 3.1).

On the other hand, Woodin (unpublished)
has shown that if the theory of $L({\Bbb
R})$ with real parameters from $V$
is frozen with respect to all set generic extensions whatsoever
then in fact
$AD^{L({\Bbb R})}$ holds, providing one more bit of evidence for the
naturalness of $AD$, the Axiom of Determinacy.

So in the light of this one obvious question is: can Woodin's  
result be strengthened by restricting the forcings to reasonable
ones?
This is non-trivial, as Woodin's proof uses the forcing for
making a singular cardinal countable, being anything but
reasonable.
More specifically: 

\bigskip
{\bf Question 1.} Suppose that for every formula $\Phi(v)$, for every
real $r \in {\Bbb R}^V$, and for every $G$ being $P$-generic where $P
\in V$ is reasonable, 
$$L({\Bbb R}^V) \models \Phi(r) {\rm \ \ iff \ \ }
L({\Bbb R}^{V[G]}) \models \Phi(r) {\rm . }$$ Does $AD^{L({\Bbb R})}$
hold (in every generic extension)?

\bigskip
We conjecture but cannot prove that the answer to this question is
"yes." In this paper, we can only give partial evidence 
in favor of this
conjecture. In fact, the argument given below can easily be
transformed to show that under its assuption global \boldmath
$\Pi$\unboldmath$^1_1$-determinacy holds.

If in particular $\omega_1$ is not to be collapsed, 
any attempt to answer this question in the affirmative seems to
essentially have to
use some coding techniques. As the coding is supposed to be set
sized, one cannot use Jensen coding as in \cite{BJW} (although it is
reasonable). A
set sized variant of it (as in \cite{ShSt}, say)
only works below $0^\sharp$ (or if $V$ is not closed under sharps, for
that matter). In general, "coding into $K$" 
techniques are called for,
where $K$ is the core model, and this is what makes the problem really
interesting (and difficult, once we get higher up in the large
cardinal zoo). 

So the above question naturally leads to the following:

\bigskip
{\bf Question 2.} Suppose that
$AD^{L({\Bbb R})}$ does not hold. Let $X \subset \omega_1$. Can $X$ be
coded into $L({\Bbb R})$ by a set-sized
reasonable forcing?

\bigskip
This would be a dual fact to the Anti-Coding Theorem of
\cite{NZ1} and \cite{NZ2}. 
However, the transit
via inner model theory for attacking this latter question
is blocked at the time of writing
by some pretty technical obstacles.   

Let us now state the main results of the present paper.

\begin{thm}\label{coding} 
Suppose that there is no inner model with a strong cardinal. 
Let $X
\subset \omega_1$. Then $X$ is \boldmath 
$\Delta$\unboldmath$^1_3$ (in
the codes) in a generic extension by a set-sized reasonable forcing.
\end{thm}

\begin{thm}\label{well-ordering}
Suppose that there is no inner model with a strong cardinal. 
Then there is a 
generic extension by a set-sized reasonable forcing
with a \boldmath
$\Delta$\unboldmath$^1_3$-well-ordering of its reals.
\end{thm}

\begin{thm}\label{freeze}
Suppose that for every $\Sigma^1_4$-formula $\Phi(v)$, 
for every $G_1$ being $P_1$-generic over $V$ where $P_1 \in V$ is
reasonable, 
for every real $r \in {\Bbb R}^{V[G_1]}$, and for every further
$G_2$
being $P_2$-generic where $P_2 \in V[G_1]$ is reasonable, 
$$V[G_1] \models \Phi(r) {\rm \ \
iff \ \ } V[G_1][G_2] \models \Phi(r) {\rm . }$$ 
Then there is an inner model
with a strong cardinal.
\end{thm}

By the following result of Woodin (unpublished),
\ref{freeze} is best possible in the sense that
one cannot derive more large cardinal strength from 
its assumption. Let $\kappa$ be a strong cardinal, 
let $\lambda = 2^{2^\kappa}$, and let $G$ be 
$Col(\lambda,\omega)$-generic over $V$.
Then any set generic extension $V[G][H_1]$ of $V[G]$ 
is $\Sigma^1_4$-correct in any of its set
generic extensions, $V[G][H_1][H_2]$. As \ref{freeze} will
be an immediate corollary to \ref{well-ordering}, \ref{well-ordering}
itself is best possible in the sense that its anti-large cardinal 
assumption cannot be weakened.

Also, suppose $V$ to be the minimal inner model with one strong
cardinal, say $V \models$ "$\kappa$ is strong," and let $G$ be 
$Col(\kappa^{++},\omega)$-generic over $V$.
By Woodin's result,
$V[G]$ is $\Sigma^1_4$-correct in its set generic extensions.
Writing $\omega_1 = \omega_1^{V[G]}$,
we have that 
${\cal J}^V_{\omega_1}$ 
cannot be \boldmath $\Sigma$\unboldmath$^1_3$ in the codes: 
otherwise the \boldmath $\Pi$\unboldmath$^1_4$-statement 
$$\forall {\rm \ countable \ } \xi \ \exists {\rm \
countable \ } \alpha > \xi \ \ \rho_\omega({\cal J}^V_\alpha) 
\leq \kappa^{++}$$
were true in $V[G]$, but could be made false by
collapsing $\omega_1$, contradicting $\Sigma^1_4$-correctness.
This implies that \ref{coding} is
best possible again
in the sense that its anti-large cardinal assumption
cannot be weakened.

Of the many problems deriving from the above questions let me
mention just one: it seems to be open whether "reasonable" can
be replaced by "proper" in the above theorems. 

\section{Coding below one strong cardinal.}

Instead of directly aiming at proving \ref{coding},
\ref{well-ordering}, and \ref{freeze} we shall first present a 
reasonable
generic extension of $V$ under the assumption that there is no inner
model with a strong cardinal. This extension will be called $V_4$
below, and it will be the case that in $V_4$ there is a
real $a$
such that
$H_{\omega_2} = J^{K(a)}_{\omega_2}$. (Here, $H_{\omega_2}$ is
the set of all sets hereditarily smaller than $\aleph_2$, and $K(a)$
is the core model built over $a$, cf. the second
next paragraph.) 
We shall then see that this
construction in fact easily gives rise to proofs of \ref{coding},
\ref{well-ordering}, and \ref{freeze}.

So let us assume thruout this section
that there is no inner model
with a strong cardinal. Then $K$, the core model below a strong
cardinal, exists (cf. \cite{J}; we here have to assume just a
little
familiarity with $K$). Moreover, $K$ is also the core
model in the sense of all set generic extensions.
We shall code an initial segment of
$V$ "into $K$." The heart of the matter will be the 
task of checking that a certain "$K$-reshaping" is
$\omega$-distributive.

For this in turn it seems to be necessary at some point during the
construction to switch from $K[A]$ to $K(A)$ for a set $A$
of ordinals, a
distinction which should be explicitly explained.
Fix $A$, a set of ordinals. Let $E$ code $K$'s extender sequence,
i.e., $K = L[E]$. Then by $K[A]$ we mean $L[E,A]$, i.e., the
constructible universe built with the two additional predicates $v \in
E$ and $v \in A$ at hand. 
Hence $K[A]$ is just
the least inner model $W$ with $K
\cup \{ A \} \subset W$. We shall also write $K_{\kappa}[A]$ for 
$J_{\kappa}[E,A]$, for an ordinal $\kappa$. Notice that the
presence of $A$ in $K[A]$ in
general destroys the internal structure of $K$. On the other hand, by
$K(A)$ we mean the core model built {\it over} 
$A$, i.e., starting from $TC(\{ A \})$,
the transitive closure of $\{ A \}$, we run
the recursive construction of $K$, with "strong $TC(\{ A \})$-mice"
instead of "strong mice" (cf. \cite{J} on the recursive definition
of $K$; to give a reference for $K(A)$, cf. 
\cite{CMIP} p. 59). By our assumption
that there is no inner model with a strong cardinal, 
$K(A)$ exists, and
in contrast to $K[A]$ it has a fine structure and can be iterated
"above $A$."

With these things in mind, we may now commence 
with our construction.       

To get things started, we use almost disjoint forcing in its
simplest form.
Fix $\delta$, a singular cardinal of uncountable cofinality and such
that $\delta^{\aleph_0} = \delta$ (for example, let $\delta$ be a
strong limit). We also may and shall assume that $\delta$ is a
cutpoint of $K$, i.e., if $E_\alpha \not= \emptyset$ is an extender
from $K$'s extender sequence with $\alpha \geq \delta$ then in
fact the citical point of $E_\alpha$ is $\geq \delta$, too. (Here we
use that there is no inner model with a strong cardinal. If there were
no such $\delta$ then using Fodor we would get a strong cardinal in
$K$.)
  
By \cite{J}, we know that $\delta^{+K} =
\delta^+$. We may also 
assume w.l.o.g. that $2^\delta = \delta^+$, because
otherwise we may collapse $2^\delta$ onto $\delta^+$ by a
$\delta$-closed preliminary forcing. We may hence
pick $A \subset \delta^+$ with the property
that $H_{\delta^+} = L_{\delta^+}[A]$. 

Now let $G_1$ be $Col(\delta,\omega_1)$-generic over $V$. Notice
that the
forcing is $\omega$-closed. 
Set $V_1 = V[G_1]$. We have that
$\omega_2^{V_1} = \delta^+ = \delta^{+K}$. Let $B \subset \omega_1$
code $G_1$ (in the sense that $G_1 \in L_{\omega_2^{V_1}}[B]$).

\bigskip
{\bf Claim 1.} In $V_1$, $H_{\omega_2} = L_{\omega_2}[A,B]$.

\bigskip
{\sc Proof.} Easy, using the fact that $Col(\delta,\omega_1)$ is
$\delta^+$-c.c.. We shall have to repeat the argument a couple of
times, so we'll be more explicit next time. 

\bigskip
\hfill $\square$ (Claim 1)

\bigskip
In what follows we let $\omega_2$ denote
$\omega_2^{V_1}$. It will also be the $\omega_2$ of all further
extensions.

Now in $K$ we may pick $(A'_\xi \ \colon \ \xi < \delta^+)$, a
sequence of almost disjoint subsets of $\delta$. In
$L_{\omega_2}[B]$ we may pick a bijective $g \colon \omega_1
\rightarrow \delta$. Then if we let $\alpha \in A_\xi$ iff $g(\alpha)
\in A'_\xi$ for $\alpha < \omega_1$ and
$\xi < \delta^+$, we have that $(A_\xi \ 
\colon \ \xi < \delta^+)$ is a sequence of almost distoint subsets of
$\omega_1$. 

In $V_1$, we may pick $A_1 \subset \omega_2$ with
$H_{\omega_2} = L_{\omega_2}[A,B] = L_{\omega_2}[A_1]$
(for example, the "join" of $A$ and $B$). We let $P_2$
be the forcing for coding $A_1$ by a subset of $\omega_1$, using the
almost disjoint sets $A_\xi$.

To be specific, $P_2$ consists of pairs $p = (l(p),r(p))$ where $l(p)
\colon \alpha \rightarrow 2$ for some $\alpha < \omega_1$ and $r(p)$
is a countable subset of $\omega_2$. We have $p = (l(p),r(p))
\leq_{P_2} q = (l(q),r(q))$ iff $l(p) \supset l(q)$, $r(p) \supset
r(q)$, and for all $\xi \in r(q)$, if $\xi \in A_1$ then $$\{ \beta
\in dom(l(p)) \setminus dom(l(q)) \ \colon \ l(p)(\beta) = 1 \} \cap
A_\xi \ = \ \emptyset.$$

By a $\Delta$-system argument, $P_2$ has the $\omega_2$-c.c. It is
clearly $\omega$-closed, so no cardinals are collapsed.
Moreover, if $G_2$ is $P_2$-generic over $V_1$, and if we set
$$C  \ = \ \bigcup_{ p \in G_2 } \ \{ \beta \in dom(l(p)) \ \colon \
l(p)(\beta) = 1 \} {\rm , }$$ then $C \subset \omega_1$
and we have that for all
$\xi < \omega_2$, $$\xi \in A_1 {\rm \ iff \ } 
Card(C \cap A_\xi) \leq
\aleph_0.$$

This means that $A_1$ is an element of any inner model containing
$(A_\xi \ \colon \ \xi < \omega_2)$ and $C$. (Of course, much more
holds.) An example of such a model is $K[C]$ in the sense explained
above. Set $V_2 = V_1[G_2]$. 
We then also have, by the same argument as for Claim 1:

\bigskip
{\bf Claim 2.} In $V_2$, $H_{\omega_2} = K_{\omega_2}[C]$

\bigskip
{\sc Proof.} Let $X$ be a subset of some 
$\gamma < \omega_2$ in $V_2$. As
$P_2$ is $\omega_2$-c.c. and $P_2 \in H_{\omega_2}$, 
$X$ has a name ${\dot X}$ in
$H_{\omega_2}^{V_1} = L_{\omega_2}[A_1]$. I.e., there is $\theta <
\omega_2$ such that ${\dot X} \in L_\theta[A_1 \cap \theta]$.
Also, there is some $\theta' < \omega_2$, $\theta' > \theta$, 
such that $(A_\xi \ \colon \
\xi < \theta) \in J_{\theta'}^K \subset 
K_{\theta'}[C]$. But then clearly ${\dot X} \in K_{\theta'}[C]$,
i.e., $X \in
K_{\omega_2}[C]$. 

This then gives $H_{\omega_2} 
\subset K_{\omega_2}[C]$ in
$V_2$. But $K_{\omega_2}[C] \subset H_{\omega_2}$ is trivial.

\bigskip
\hfill $\square$ (Claim 2)

\bigskip
Now let $C_1 \subset \omega_1$ 
be such that $C \in L_{\omega_2}[C_1]$ as well as ${\cal
J}^K_\delta \in L_{\omega_2}[C_1]$ (for example, the "join" of $C$ and
a code for ${\cal J}^K_\delta$).
     
Our task is now to code "down to a real," i.e.,
we want to find a further ($\omega$-distributive) 
generic extension in which
$H_{\omega_2} \subset K(a)$ where $a$ is a real. As we cannot expect
$\omega_1$ to be a successor in $K$ (and as we certainly cannot force
this 
to be the case using an $\omega_1$-preserving
forcing), we have to use the slightly
more advanced coding technique which first requires $C_1$ to become
"reshaped." 

However, there is no hope of showing that reshaping is
reasonable
if we work with the wrong definition of "being reshaped."
Thus, the reader should notice the round brackets in the following
definition. 

\begin{defn}
Let $X \subset \omega_1$. We say that an $f$ is $X$-reshaping if 
$f \colon \alpha \rightarrow 2$ for some $\alpha \leq \omega_1$ and
moreover for all $\beta \leq
\alpha$, $$K(X \cap \beta , f \upharpoonright \beta) \models
Card(\beta) \leq \aleph_0.$$ 
\end{defn}

Now let $P_3$ be the forcing for adding a $C_1$-reshaping 
characteristic function of a subset of 
$\omega_1$. Formally, $p \in P_3$ iff $p$ is $C_1$-reshaping and
$dom(p) < \omega_1$. The order is by reverse inclusion, i.e., $q
\leq_{P_3} p$ iff $q \supset p$. 
(In fact we may assume "w.l.o.g." that every
$C_1$-reshaping $p$ has $dom(p) < \omega_1$, because otherwise we
could just fix a counterexample and go ahead with
forcing with $P_4$.)

It is easy to see that for any $\alpha < \omega_1$, the set $D^\alpha
= \{ p \in P_3 \ \colon \ dom(p) \geq \alpha \}$ is open dense in
$P_3$. In fact, given $q \in P_3$ with $dom(q) < \alpha$,
we may let $p \leq q$ with $dom(p) = \alpha + \omega$ 
be such that $p
\upharpoonright [dom(q),dom(q) + \omega)$ codes some 
bijective $f \colon \omega
\rightarrow \alpha + \omega$, and then $p \in D^\alpha$.  

\bigskip
{\bf Claim 3.} $P_3$ is $\omega$-distributive.

\bigskip
{\sc Proof.} We first
need the following observation, which is due to John Steel.
The argument proving it also has to be repeated a couple of times.         

\bigskip
{\bf Subclaim 1.} (Steel)
$H_{\omega_2} = J^{K(C_1)}_{\omega_2}$.

\bigskip
{\sc Proof.} Let $W = K^{K(C_1)}$, the core model built inside 
$K(C_1)$.
In $V_2$, 
let $\theta$ be any regular cardinal, and let $\kappa$ be a
singular cardinal with $cf(\kappa) > \theta$. 
By weak
covering applied inside $K(C_1)$ (cf. \cite{J}), we have that
$cf^{K(C_1)}(\kappa^{+W}) \geq \kappa$, which implies that
$cf^{V_2}(\kappa^{+W}) \geq cf^{V_2}(\kappa) > \theta$.
Hence, inside $V_2$, 
for any regular $\theta$ there is a stationary class of
cardinals $\kappa$ with $cf(\kappa^{+W}) > \theta$. 
  
But this implies that inside $V_2$, 
$W$ is a universal weasel in the sense 
that the coiteration of $W$ with any (set sized) 
mouse terminates after
$< OR$ many steps and the mouse-side is simple. 

Moreover, as ${\cal
J}^K_\delta \in K(C_1)$, an easy absoluteness argument using the
recursive definition of $K$ yields that ${\cal J}^K_\delta
\triangleleft W$, i.e., ${\cal J}^K_\delta$ is an initial segment of
$W$. Let us consider the coiteration of $K = K^{V_2}$ with $W$.
 
Because $\delta$ is a cutpoint in $K$ and ${\cal J}^K_\delta
\triangleleft W$, 
the coiteration is above $\delta$ on the $K$-side. In fact, it must be
above $\delta^{+K} = \omega_2$ on the $K$-side, as otherwise we would
have to "drop" to a mouse on the $K$-side and by non-soundness
of the further iterates of $K$
the coiteration would have to last
$OR$ many steps, contradicting the universality of $W$.

But the coiteration has to be above $\delta$ on the $W$-side, too,
because otherwise
we may replace $K$ by a very soundness witness for a 
large enough initial segment of $K$ and use its definability
property
everywhere below $\delta$ to get a contradiction as usual. But then 
the coiteration must be above $\delta^{+W}$ on the $W$-side, 
too, this time
by the universality of $K$.

This now means that ${\cal
J}^W_{\delta^{+W}} = {\cal J}^K_{\delta^{+K}}$, and hence
$J^K_{\omega_2} \subset K(C_1)$. 
In particular, $(A_\xi \ \colon \ \xi
< \omega_2) \in K(C)$, in fact $(A_\xi \ \colon \ \xi < \gamma)
\in J^{K(C_1)}_{\omega_2}$ for any $\gamma < \omega_2$ by
acceptability, which easily gives the claim.

\bigskip
\hfill $\square$ (Subclaim 1)

\bigskip
We remark in passing that Subclaim 1 would not have to hold
if $\delta$ had not been chosen as a cutpoint of $K$.

We now fix a condition $p \in P_3$, and open dense sets
$D'_i$, $i < \omega$. We have to find $q \leq_{P_3} p$ 
with $q \in D'_i$ for
every $i < \omega$.  

We may pick $$\pi \ \colon \ N \rightarrow_{\Sigma_1} 
{\cal J}^{K(C_1)}_{\omega_2} =
H_{\omega_2} {\rm , }$$
where $N$ is countable, $\kappa = c.p.(\pi)$ is such that
$\rho_1(N) = \kappa$, 
$N$ is sound above $\kappa$, and
$\{ p \} \cup \{ D'_i \ \colon \ i < \omega \} \subset ran(\pi)$. 
Notice that $\pi(\kappa) = \omega_1$.

This situation is obtained, for example, if we first let
$M$ be
the $\Sigma_1$-hull of $\omega_1 \cup \{ p \} \cup \{ 
(D'_i \ \colon \ i < \omega) \}$, taken
inside ${\cal J}^{K(C_1)}_{\omega_2}$, and then let
$\pi$ be ${\tilde \pi} \upharpoonright N$ for some
${\tilde \pi} \ \colon \ {\tilde N} \rightarrow H_{\omega_2}$ with $M
\in ran({\tilde \pi})$, $N = {\tilde \pi}^{-1}(M)$,
and ${\tilde N}$ being countable. 

\bigskip
{\bf Subclaim 2.} $N \in K(C_1 \cap \kappa)$.

\bigskip
{\sc Proof.} We coiterate $N$ with $K(C_1 \cap \kappa)$, getting
comparable $N^*$ and $K^*$. As $\rho_1(N) = \kappa$, every
non-trivial iterate of $N$ is non-sound. Hence, by the universality
of $K(C_1 \cap \kappa)$, $N$ cannot be moved at all in the comparison,
i.e., $N^* = N$. 
  
Now suppose that $K(C_1 \cap \kappa)$ were to be moved. 
If there is a
drop, then $K^*$ is non-sound, 
and hence $N \triangleleft K^*$ by the soundness of $N$. But
if the iteration is simple then we trivially have $N \triangleleft
K^*$ as well. Because we assume $K^* \not= K(C_1 \cap \kappa)$, 
letting $\nu$ be the index
of the first extender of the iteration from $K(C_1 \cap \kappa)$ to
$K^*$, we have that $\nu \leq OR \cap N$ and $\nu$ is a cardinal in
$K^*$. But this is a contradiction, as $N \triangleleft K^*$ and
$\rho_1(N) = \kappa$, so that
$K^*$ knows that there are no cardinals in the half-open interval
$(\kappa,OR \cap N]$.

This shows that in fact $N \triangleleft K(C_1 \cap \kappa)$, which in
particular gives us what we want.

\bigskip
\hfill $\square$ (Subclaim 2)  

\bigskip

Now let $\{ D_i \ \colon \ i < \omega \}$ be the open dense sets in
$ran(\pi)$. It suffices to construct $q \leq_{P_3} 
p$ with $q \in D_i$ 
for every $i < \omega$. For this we use an argument of \cite{ShSt}. 

We may assume w.l.o.g. that $\kappa = \omega_1^{K(C_1 \cap
\kappa)}$, as otherwise $\kappa$ is countable in $K(C_1 \cap \kappa)$
and the task of constructing $q$ turns out to be trivial.
But nevertheless $N$ has size $\kappa$ in $K(C_1 \cap \kappa)$ (because
$\rho_1(N) = \kappa$). Hence we may pick a club $E \subset \kappa$ in
$K(C_1 \cap \kappa)$ which grows faster than all clubs in $N$, i.e.,
whenever ${\bar E} \subset \kappa$ is a club in $N$ then $E \setminus
{\bar E}$ is bounded in $\kappa$. 

Inside $K(C_1 \cap \kappa)$, 
we are going to construct a sequence $(p_i \ \colon
\ i < \omega)$ of conditions below $p$
such that $p_{i+1} \leq_{P_3} p_i$ and
$p_{i+1} \in D_i$. We also want to maintain inductively
that $p_{i+1} \in N$. (Notice that $p \in N$
to begin with.) 
In the end we also want to have that setting $q =
\cup_{i < \omega} \ p_i$, we have that
$q \in P_3$, which of course is the the
non-trivial part.

To commence, let $p_0 = p$. Now suppose that $p_i$ is given, $p_i \in
N$. Set
$\alpha = dom(p_i) < \kappa$. Work inside $N$ for a minute.
For all $\beta$ such that $\alpha \leq
\beta < \kappa$ we may pick some $p^\beta
\leq_{P_3} p_i$ such
that: $p^\beta \in \pi^{-1}(D_i)$, $dom(p^\beta) > \beta$,
and for all limit ordinals
$\lambda$, $\alpha \leq \lambda \leq \beta$, $p^\beta(\lambda)
= 1$ iff $\lambda = \beta$.  
Then there is ${\bar E}$ club in $\kappa$ such that for any $\eta \in
{\bar E}$, $\beta < \eta \Rightarrow dom(p^\beta) <
\eta$. 

Now back in $K(C_1 \cap \kappa)$, we may pick 
$\beta \in E$ such that $E \setminus {\bar E} \subset
\beta$. Set $p_{i+1} = p^\beta$, and let for future reference $\beta =
\beta_{i+1}$.
Of course $p_{i+1} \in D_i \cap N$.
We also have that $dom(p_{i+1}) < min \{ \epsilon \in E \ \colon \
\epsilon > \beta \}$, so that
for all limit ordinals $\lambda \in E
\cap (dom(p_{i+1}) \setminus dom(p_i))$ we have that 
$p_{i+1}(\lambda) = 1$ iff $\lambda = \beta_{i+1}$.

Now set $q = \cup_{i < \omega} \ p_i$. 
Well, for every $\alpha < \kappa$, $D^\alpha = \{ r \in P_3 \ \colon \
dom(r) \geq \alpha \} \in ran(\pi)$, and hence $dom(q) \geq \kappa$.
Alas, we also have $dom(q) \leq \kappa$, because $p_i \in N$ and so 
$dom(p_i) < \kappa$ for all
$i < \omega$. Hence we have arranged that $dom(q) = \kappa$.

We are done if we can show that
$q$ is a condition. The only problem here is 
to show that $$K(C_1 \cap
\kappa,q) \models Card(\kappa) \leq \aleph_0.$$ 
But by the construction
of the $p_i$'s we have that $$\{ \lambda \in E \cap (dom(q) \setminus
dom(p)) \ \colon \ \lambda {\rm \ is \ a \ limit \ ordinal \ and \ }
q(\lambda) = 1 \}$$ $$= \ 
\{ \beta_{i+1} \ \colon \ i < \omega \} {\rm ,
}$$ being a cofinal subset of $E$.
It hence suffices to verify that
$E$ is an element of $K(C_1 \cap \kappa , q)$, because then
$\{ \beta_{i+1} \ \colon \ i < \omega \} \in K(C_1 \cap \kappa , q)$
witnesses that $Card(\kappa) \leq \aleph_0$.

But $E \in K(C_1 \cap \kappa , q)$ is shown by the argument for
subclaim 1: we build $W = K(C_1 \cap \kappa)^{K(C_1 \cap \kappa,q)}$,
observe its universality in $V_2$, and deduce that $E \in
J_{\omega_2}^{K(C_1 \cap \kappa)} = J_{\omega_2}^W \subset K(C_1 \cap
\kappa , q)$.

\bigskip
\hfill $\square$ (Claim 3)

\bigskip
Actually, there is an easier proof showing
$\omega$-distributivity of $P_3$ than the 
one we just gave. We however chose to
include the above argument because it gives more, namely it "almost"
proves properness of $P_3$:

\bigskip
{\bf Claim 3'.} Let $S \in {\cal P}(\omega_1) \cap V_2$ be stationary
in $V_2$. Then $P_3 \forces$ "${\hat S}$ is stationary."

\bigskip
{\sc Proof.} Let $r \in P_3$, $r \forces$ "${\dot C}$ is a club subset
of $\omega_1$." In the above
argument, let $p \leq_{P_3} r$, and $\pi \ \colon \ N \rightarrow
H_{\omega_2}$ be such that ${\dot C} \in ran(\pi)$ and $\kappa \in
S$.

But then for the $q$ constructed we get that $q \forces$
"${\hat \kappa} \in {\hat
S} \cap {\dot C}$."

\bigskip
\hfill $\square$ (Claim 3')

\bigskip
Unfortunately, I cannot decide whether $P_3$ always has to be proper.
We'll return to this issue at the end of this section.    

Now let $G_3$ be $P_3$-generic over $V_2$, and set $V_3 = V_2[G_3]$.
Let $C'$ be that subset of $\omega_1$ having $\cup_{p \in
G_3} \ p$ as its characteristic function, and let $D$ be the "join" of
$C$ and $C'$. Again we get:

\bigskip
{\bf Claim 4.} In $V_3$, $H_{\omega_2} = J^{K(D)}_{\omega_2}$.

\bigskip
{\sc Proof.} The only new point here is that we have to check that
$J_{\omega_2}^{K(C)} \subset J_{\omega_2}^{K(D)}$. But here we can
argue exactly as at the end of the proof of claim 3, 
by building $K(C)^{K(D)}$, observing
its universality, and deducing $J_{\omega_2}^{K(C)} \subset K(D)$.

\bigskip
\hfill $\square$ (Claim 4)

\bigskip   
We may now finally code down to a real by using almost disjoint
forcing once more. By the fact that $D$ is "$K$-reshaped," there is a
(unique) sequence $(a_\beta \ \colon \ \beta 
< \omega_1)$ of 
subsets of $\omega$ such that for each $\beta 
< \omega_1$, $a_\beta$ is
the $K(D \cap \beta)$-least subset of $\omega$ being almost disjoint
from any $a_{\bar \beta}$ for ${\bar \beta} < \beta$.

We then let $P_4$ consist of all pairs $p = (l(p),r(p))$ where
$l(p) \colon n \rightarrow 2$ for some $n < \omega$ and $r(p)$ is a
finite subset of $\omega_1$. We let $p = (l(p),r(p)) \leq_{P_4} q =
(l(q),r(q))$ iff $l(p) \supset l(q)$, $r(p) \supset r(q)$, and for all
$\beta \in r(q)$, if $\beta \in D$ then $$\{ \gamma \in dom(r(p))
\setminus dom(r(q)) \ \colon \ r(p)(\gamma) = 1 \} \cap a_\beta \ = \
\emptyset.$$ 

By another $\Delta$-system argument, $P_4$ has the c.c.c..
Moreover, if $G_4$ is $P_4$-generic over $V_3$, and if we set
$$a = \bigcup_{p \in G_4} \ \{ \gamma \in dom(l(p)) \ \colon \
l(p)(\gamma) = 1 \} {\rm , }$$ then we have that for $\gamma <
\omega_1$, $$\gamma \in D {\rm \ iff \ } Card(a \cap a_\gamma) <
\aleph_0.$$ 

We 
finally get: 

\bigskip
{\bf Claim 5.} In $V_4$, $H_{\omega_2} = J_{\omega_2}^{K(a)}$.

\bigskip
{\sc Proof.} In order to show that $J_{\omega_2}^{K(D)} \subset K(a)$
we first have to verify $D \in K(a)$. This this is easily seen to
follow (combined with the uniform definability of the $a_\beta$'s)
from the fact that if $D \cap \gamma \in K(a)$ then
$J_{\omega_2}^{K(D \cap \gamma)} \subset K(a)$ (which in turn 
is true by the
argument for Subclaim 1).

But then one more Subclaim 1 type argument proves $J_{\omega_2}^{K(D)}
\subset K(a)$.

\bigskip
\hfill $\square$ (Claim 5)  

\bigskip
At last, we observe that the $4$-step iteration yielding from $V$ to
$V_4$ is reasonable. Let $\alpha$ be any infinite ordinal.
Because $P_1$, $P_2$, and $P_3$ are all
$\omega$-distributive, we have that $[\alpha]^\omega \cap V = 
[\alpha]^\omega \cap V_3$. But $P_4$ is c.c.c., hence proper, which
implies that $[\alpha]^\omega \cap V$ is stationary in
$[\alpha]^\omega \cap V_4$.

Actually, we also get that every $S \in {\cal P}(\omega_1) \cap V$
which is stationary in $V$ remains stationary in $V_4$. Clearly, any
such $S$ is still stationary in $V_2$, as $P_1$ and $P_2$ are both
$\omega$-closed, hence proper. But then $S$ is stationary in $V_3$ by
Claim 3', and so in the end $S$ is still stationary in $V_4$ by the
properness of $P_4$.

\section{Getting the theorems.}

We may now easily derive \ref{coding}, \ref{well-ordering}, and
\ref{freeze} from the work done in the previous section.

\bigskip
{\sc Proof} of \ref{coding}. Let us assume that
there is no inner model with a strong cardinal. 
Fix $X \subset \omega_1$.
Running the
construction of the last section, we may certainly assume w.l.o.g.
that for all $\xi < \omega_1$, $\xi \in X$ iff $2 \xi \in D$. But then
we have that in $V_4$,
$\xi \in X$ iff $$\exists (a_\gamma \colon \gamma \leq 2
\xi) \ \exists d \subset 2 \xi + 1 \ [ \ \forall \gamma \leq 2 \xi 
\ ( \
a_\gamma {\rm \ is \ the \ } K(d \cap \gamma) {\rm -least \ }$$ 
$${\rm subset \ of \ } \omega {\rm \ a. \ d. \ from \ all \ } 
a_{\bar \gamma}
{\rm
, \ } {\bar \gamma} < \gamma {\rm , }$$  
$${\rm \ and \ } \gamma \in d \leftrightarrow
Card(a \cap a_\gamma) < \aleph_0 \ ) \ 
{\rm \ and \ } 2 \xi \in d \
].$$
But $a_\gamma$ is the $K(d \cap \gamma)$-least subset of $\omega$
a. d. from all $a_{\bar \gamma}$, ${\bar \gamma} < \gamma$, iff
$a_\gamma$ is the ${\cal M}$-least subset of $\omega$ a. d. from all
$a_{\bar \gamma}$, ${\bar \gamma} < \gamma$, for all $(d \cap
\gamma)$-mice ${\cal M}$ with $\{ a_{\bar \gamma} \colon {\bar \gamma}
\leq \gamma \} \subset {\cal M}$. 
We may hence rewrite the displayed formula as
$$\exists (a_\gamma \colon \gamma \leq 2 \xi) \
\exists d \subset 2 \xi + 1 \ \exists ({\cal M}_\gamma \colon 
\gamma \leq 2 \xi) \ [ \ \forall \gamma \leq 2 \xi \
( \ {\cal M}_\gamma {\rm \ is \ a }$$ 
$$(d \cap \gamma) {\rm -mouse \ with \ } 
\{ a_{\bar \gamma} \colon {\bar \gamma} \leq \gamma \} \subset
{\cal M}_\gamma {\rm , }$$
$$a_\gamma {\rm \ is \ the \ }
{\cal M}_\gamma {\rm -least \ subset \ of \ } \omega {\rm \ a. \ d. \
from \ all \ } a_{\bar \gamma} {\rm , \ } {\bar \gamma} < \gamma 
{\rm , }$$ $${\rm and \ } \gamma \in d \leftrightarrow Card(a \cap
a_\gamma) < \aleph_0 \ ) \ {\rm and \ } 2 \xi \in d \ ].$$ 
As "mousehood" is $\Pi^1_2$, this
latter displayed formula is now easily be seen to be $\Sigma^1_3(a)$.  
Hence $X$ turns out to be $\Sigma^1_3(a)$ in
the codes, inside $V_4$.

But because $(a_\gamma \colon \gamma < \omega_1)$ is uniquely
determined given $a$, it is easy that $X$ is now $\Pi^1_3(a)$ 
in the codes as well, inside $V_4$.

\bigskip
\hfill $\square$ (\ref{coding})

\bigskip
{\sc Proof} of \ref{well-ordering}. Assuming that there is no inner
model with a strong cardinal and building $V_4$ as in the last
section, we have by Claim 5 that $H_{\omega_2} = J_{\omega_2}^{K(a)}$
in $V_4$. Inside $V_4$, we may thus define a well-ordering of the
reals by $x < y$ iff $$\exists {\cal M} \ [ \ {\cal M} {\rm \ is \ 
an \ } a {\rm -mouse \
and \ } x <_{\cal M} y \ ].$$ (Here, $<_{\cal M}$ denotes the
order of constructibility of ${\cal M}$.) Again using the fact that
"mousehood" is $\Pi^1_2$, this gives a $\Sigma^1_3(a)$-wellordering of
${\Bbb R}$. But clearly $x < y$ can now also written in a
$\Pi^1_3(a)$-fashion.

\bigskip
\hfill $\square$ (\ref{well-ordering})

\bigskip
{\sc Proof} of \ref{freeze}. Let us first suppose that there is a
\boldmath $\Pi$\unboldmath$^1_3$-well-ordering of ${\Bbb R}$. This
fact can be expressed by a \boldmath
$\Pi$\unboldmath$^1_4$-statement $\Phi$
being true in $V$.
By
adding $\omega_1$ many Cohen reals with finite support
we get a (c.c.c., hence
proper, hence) reasonable extension
$V[G]$ of $V$ in which there is no projective well-ordering of the
reals, so that $\Phi$ fails.

Thus if we assumption of \ref{freeze} holds then there can be no
\boldmath $\Pi$\unboldmath$^1_3$-well-ordering of the reals in $V$.
Now let us suppose that there is no inner model with a strong
cardinal. 

Then by \ref{well-ordering} there is a reasonable extension $V[G]$
of $V$
with a \boldmath $\Pi$\unboldmath$^1_3$- (in fact
\boldmath $\Delta$\unboldmath$^1_3$-) well-ordering of its reals.
As above, this fact can be expressed by a
\boldmath $\Pi$\unboldmath$^1_4$-statement
$\Phi$ being true in $V[G]$. But again by adding $\omega_1$ many Cohen
reals we get a reasonable extension $V[G][H]$ of $V[G]$ in which there
is no projective well-ordering of the reals, so that $\Phi$ fails.
Contradiction!

We have shown that there is an inner model with a strong cardinal.

\bigskip
\hfill $\square$ (\ref{freeze}) 

\section{Further remarks.}

In this last section we want to point out a couple of things which are
of related interest. The first one is an immediate corollary to
\ref{well-ordering}, if combined with an old argument of Harrington
(cf. \cite{jech} pp. 575 ff.).

\bigskip
{\bf Fact 1.}
Suppose that there is no inner model with a strong cardinal. Let
$\delta$ be any ordinal. Then there is a generic extension by a
set-sized reasonable forcing in which \boldmath
$\delta$\unboldmath$^1_4 > max \{ \delta , \aleph_3 \}$.

\bigskip
If $V$ is not closed under $\sharp$'s then a slight variation of the
Harrington argument 
actually shows that we can make 
\boldmath $\delta$\unboldmath$^1_3$ as
large as prescribed by a set-sized 
reasonable forcing. On the other hand, if
${\Bbb R}$
is closed under $\sharp$'s (i.e., if \boldmath
$\Pi$\unboldmath$^1_1$-determinacy holds)
then \boldmath $\delta$\unboldmath$^1_3
\leq \aleph_3$.

\bigskip
{\sc Proof} of Fact 1 (sketch).
Going first to a reasonable extension provided by
\ref{well-ordering}, we then run the Harrington construction as in
\cite{jech}, pp. 575 ff., say. 
As the two further forcings are c.c.c., hence proper,
the final
extension is a reasonable one. 

\bigskip
\hfill $\square$ (Fact 1, sketch)
  
\bigskip
Finally, we want to point out a "converse" to \ref{well-ordering}.

\bigskip
{\bf Fact 2.}
Suppose that there is no inner model with a strong cardinal, and
${\Bbb R}$ is closed under $\sharp$'s (i.e., \boldmath
$\Pi$\unboldmath$^1_1$-determinacy holds). Let $a \in {\Bbb R}$.
If there is a
$\Sigma^1_3(a)$-well-ordering of ${\Bbb R}$ then in fact
${\Bbb R} = {\Bbb R} \cap K(a)$.  
                 
\bigskip
This
is false if ${\Bbb R}$ is not closed under $\sharp$'s. On the
other hand, we can weaken the anti-large cardinal assumption to "there
is no inner model with a Woodin cardinal" if in addition we are
willing to assume the existence of two measurable cardinals (in $V$)
(this follows from Steel's $\Sigma^1_3$-correctness theorem, 
cf. \cite{CMIP}).

\bigskip
{\sc Proof} of Fact 2 (sketch). 
Under the hypotheses of the theorem, by \cite{StW} there
is a tree $T_2 \in K(a)$ projecting to the universal \boldmath
$\Pi$\unboldmath$^1_2$-set of reals. (This tree gives
$\Sigma^1_3$-correctness of $K(a)$ in $V$.)

But we can now run the Mansfield proof that if there is a
$\Sigma^1_2(b)$-well-ordering of the reals then in fact ${\Bbb R}
\subset L[b]$, for $b \in {\Bbb R}$ (cf. \cite{jech} Theorem 100
(c), p. 534), but with the
Shoenfield tree replaced by $T_2$. This gives the desired result.

\bigskip
\hfill $\square$ (Fact 2, sketch)

\bigskip
As a by-product of this proof we get that if $A$ is a
$\Sigma^1_3(a)$-set of reals, $a \in {\Bbb R}$, and $A \setminus K(a)
\not= \emptyset$, then $A$ in fact has a perfect subset. Hence if
$\omega_1^{K(a)}$ is countable for every $a \in {\Bbb R}$ then every
\boldmath $\Sigma$\unboldmath$^1_3$-set of reals has the perfect
subset property.

{\sc Mathematisches Intitut, Uni Bonn, Beringstr. 4, 53115 Bonn,
Germany,} 

{\tt rds@math.uni-bonn.de}

\bigskip
{\sc Math Department, UCB, Berkeley CA 94720, USA,} 

{\tt
rds@math.berkeley.edu}

\end{document}